\documentclass[12pt]{article}
\usepackage{amsmath}
\usepackage{amscd}
\usepackage{amssymb}
\usepackage{array}

\newtheorem{definition}{Definition}[section]
\newtheorem{theorem}[definition]{Theorem}

\newtheorem{proposition}[definition]{Proposition}

\newtheorem{lemma}[definition]{Lemma}

\begin{document}

\rightline{IFIC/04-24}

\rightline{FTUV-04/0521}

\newcommand{\R}{\mathbb{R}}
\newcommand{\C}{\mathbb{C}}
\newcommand{\Z}{\mathbb{Z}}
\newcommand{\N}{\mathbb{N}}
\newcommand{\bA}{\mathbb{A}}

\newcommand{\g}{\mathcal{G}}
\newcommand{\A}{\mathcal{A}}
\newcommand{\B}{\mathcal{B}}
\newcommand{\PP}{\mathcal{P}}
\newcommand{\I}{\mathcal{I}}
\newcommand{\La}{\mathcal{L}}

\newcommand{\cFA}{\mathcal{FA}}

\newcommand{\m}{\mathbf{m}}

\newcommand{\de}{\delta}

\newcommand{\cR}{\mathcal{R}}

\newcommand{\defi}{\scriptsize{\mbox{def}}}
\newcommand{\CX}{C^{\infty}(X)}
\newcommand{\CAn}{\C[\bA^n]}

\newcommand{\Ah}{\A_h}
\newcommand{\dex}{\partial_x}
\newcommand{\dey}{\partial_y}

\newcommand{\lra}{\longrightarrow}

\newcommand{\Span}{\rm{span}}
\newcommand{\Sym}{\rm{Sym}}
\newcommand{\LT}{\rm{LT}}
\newcommand{\rmod}{\rm{mod}}
\newcommand{\rId}{\mathrm{Id}}

\vskip 2cm \centerline{\LARGE \bf On the deformation quantization}
\bigskip
\centerline {\LARGE\bf of affine   algebraic varieties}

\vskip 1cm

\centerline{R. Fioresi$^\flat$\footnote{Investigation supported by
the University of Bologna, funds for selected research topics.},
M. A. Lled\'o$^\natural$ and V. S. Varadarajan $^\sharp$}

\bigskip

\centerline{\it $^\flat$ Dipartimento di Matematica, Universit\`a
di Bologna }
 \centerline{\it Piazza di Porta S. Donato, 5. 40127 Bologna. Italy.}
\centerline{{\footnotesize e-mail: fioresi@dm.unibo.it}}

\bigskip

\centerline{\it $^\natural$ Departament de F\'{\i}sica Te\`orica,
Universitat de Val\`encia. }
 \centerline{\it C/ Dr. Moliner, 50. 46100 Burjassot (Val\`encia) Spain.}
\centerline{{\footnotesize e-mail: maria.lledo@ific.uv.es}}

\bigskip

\centerline{\it $^\sharp$ Department of Mathematics, UCLA.}
\centerline{\it Los Angeles, CA, 90095-1555, USA}
\centerline{{\footnotesize e-mail: vsv@math.ucla.edu}}

\vskip 2cm
\begin{abstract} We compute an explicit algebraic deformation
quantization for an affine Poisson variety described by an ideal
in a polynomial ring, and inheriting its Poisson structure from
the ambient space.

\end{abstract}

\vfill\eject

\section{ Introduction }

Since the fundamental work of Bayen et al \cite{bffls} in the
seventies, a lot of effort has been dedicated to show the
existence of deformations of a Poisson manifold. Some landmarks in
this way were the proof of the existence of differential star
products for symplectic manifolds which was done independently by
De Wilde and Lecomte \cite{dl} and Fedosov \cite{fe}, using
different constructions. It turned out that the star products on a
symplectic manifold are classified, up to equivalence, by the de
Rahm cohomology $H^2(M)$. Etingof and Kazhdan showed the existence
of star products for another class of Poisson manifolds, the
Poisson-Lie groups. Kontsevich gave the proof of existence and
classification of star products on an arbitrary Poisson manifolds
as a consequence of his formality theorem \cite{ko1}. Tamarkin
\cite{ta} gave another proof of the formality theorem that relates
it to Deligne's conjecture on Hochschild's complexes.

More recently, there has grown an interest in translating all the
results mentioned above (valid for $C^\infty$ manifolds and
differential deformations) to the algebraic geometric setting
\cite{ko2,ye,bk}. We will comment on these approaches in Section
\ref{preliminaries}. On the other hand, algebraic star products on
the sphere appeared as soon as in Refs. \cite{cg}, and later on, a
more general construction appeared in Refs. \cite{fl,ll}. In these
last references, the importance in physics of algebraic (and not
necessarily differential) star products  was stressed, because
they are the physical choice in problems as fundamental as the
quantization of angular momentum. As a consequence, they are also
related to geometric quantization. This was not taken into account
the original papers when the differentiability hypothesis was
assumed through the whole process of deformation quantization.

The relation between algebraic and differential star products was
intriguing, and it was studied in successive papers
\cite{fll,fl2}. The approach followed in these articles was
restricted to coadjoint orbits of semisimple groups (so, to linear
Poisson structures). Explicit algebras, defined by generators and
relations, where considered and this  allowed  to find some new
features in the algebraic case. One could easily see that there
were non isomorphic deformations even in the simplest case.

\bigskip

It is our intention in this work to extend the approach of of
Refs.  \cite{fl,ll} to a wider case of affine Poisson algebraic
varieties, where the degree of the Poisson structure is arbitrary.
Although working with a restricted class of algebraic varieties,
the main advantage of our approach is that the deformation is
shown explicitly in terms of generators and relations, with no
recursion to gluing procedures. We will construct a suitable non
commutative algebra and then we will show that it is an algebraic
deformation of the coordinate ring of the variety. The precise
sense of this statement is explained in Section
\ref{preliminaries}, where we fix the notation and briefly explore
other approaches present in the literature. The starting point is
the deformation of the affine space given by Kontsevich
\cite{ko1}. In Section \ref{affinespace} we present such algebra
as a quotient of the tensor algebra by a certain ideal, very much
in the way that the universal enveloping algebra is presented
(which is the deformation of a linear Poisson structure
\cite{ho}). In the case of the coadjoint orbits, the deformed
algebra was obtained by quotienting the enveloping algebra by an
ideal that it is related to the ideal of the variety. In Section
\ref{algvariety} the same procedure is extended to a bigger class
of algebraic varieties, with no particular restriction on the
degree of the Poisson structure. Some assumptions, nevertheless,
must be made in the course of the proof, but they are of technical
nature and it is likely that they can be dropped. At this  moment
we do not know if this is possible.

Furthermore, one may think on gluing the deformations obtained for
affine varieties to deform more general algebraic varieties,
perhaps with a procedure \`a la Fedosov. This is out of the scope
of the present paper but may be approached in other works.


\section{ Preliminaries \label {preliminaries}}

In this section we want to introduce some of the key definitions
of the theory of deformation quantization. In particular we want
 to compare our definitions and approach with the ones appearing
 in the literature.

\begin{definition}\label{formaldef}
Let $(\A,\{\;,\,\})$ be a Poisson algebra over a field $k$. We say
that the associative algebra $\A_{[h]}$ over $k[[h]]$ is a formal
deformation of $\A$ if
\medskip

\noindent 1. There exists an isomorphism of $k[[h]]$-modules
$\psi: \A[[h]]\longrightarrow \A_{h}$;

\noindent 2.  $\psi(f_1f_2)=\psi(f_1)\psi(f_2)\; {\rm
mod(h)},\quad \forall f_1, f_2 \in A[[h]]$;

\noindent 3.
$\psi(f_1)\psi(f_2)-\psi(f_2)\psi(f_1)=h\psi(\{f_1,f_2\})\;\; {\rm
mod}(h^2), \quad \forall f_1,f_2\in \A[[h]]$. \hfill$\blacksquare$
\end{definition}

If $\A_{\C}$ is the complexification of a real Poisson algebra
$\A$ we can give the definition of {\it formal deformation} of
$\A_{\C}$ by replacing $\R$ with $\C$ in Definition
\ref{formaldef}. A reality condition on the star product may be
required in order to have a star product defined also over $\R$.

The associative product in $\A[[h]]$ defined by
\begin{equation}
f \star g= \psi^{-1}(\psi(f)\cdot\psi(g)), \qquad f,g\in
\A[[h]]\label{isom}
\end{equation}
is called the {\it star product on $\A[[h]]$ induced by $\psi$.}

A star product on $\A[[h]]$ can be also defined as an associative
$k[[h]]$-linear product given by the formula \begin{equation} f
\star g=fg+B_1(f,g)h+B_2(f,g)h^2+\dots \in \A[[h]],\quad f,g\in \A
\label{starproduct}\end{equation} where the $B_i$'s are bilinear
operators. The associativity of $\star$ implies that
$\{f,g\}=B_1(f,g)-B_1(g,f)$ is a Poisson bracket on $\A$. So this
definition is a special case of the previous one where
$\A_h=\A[[h]]$ and $\star$ is induced by $\psi={\rm Id}$.

Two star products on $\A[[h]]$, $\star$ and $\star'$  are said to
be {\it equivalent} (or {\it gauge equivalent}) if there exists a
linear map $T:\A[[h]]\rightarrow\A[[h]]$ of the form
$$T=\rId+\sum_{n
>0}h^nT_n$$ with $T_n$ linear operators on $\A[[h]]$, such that
$$
f\star g=T^{-1}(T(f)\star'T(g)).
$$
Two star products that are equivalent are isomorphic and have the
same first order term, so they are formal deformations of the same
Poisson structure.

If $\A \subset C^{\infty}(M)$ and the operators $B_i$'s are
bidifferential operators we say that the star product is {\it
differential}. If in addition $\A=C^{\infty}(M)$ and $M$ is a real
Poisson manifold, we will say that $\star$ is a {\it differential
star product on $M$}. The set of  (gauge) equivalence classes of
differential star products on a manifold $M$ has been classified
by Kontsevich
  in terms of equivalence classes of formal Poisson
 structures (modulo formal diffeomorphisms) \cite{ko1}.

 \medskip

Let $M_{\C}$ be a complex algebraic affine variety defined over
$\R$, whose real points are  a real algebraic Poisson variety $M$.
We denote by $\A_{\C}=\C[M_{\C}]$ its  coordinate ring, which is a
Poisson algebra. We will say that a formal deformation of
$\A_{\C}$ is a formal deformation of the algebraic variety
$M_{\C}$.

We have an {\it algebraic star product on $M$} (or $M_{\C}$) if
the bilinear operators $B_i$ in (\ref{starproduct}) are algebraic
operators.

A special example of algebraic Poisson varieties are the coadjoint
orbits Lie groups. They are the symplectic leaves of the Kirillov
Poisson structure in the dual of the Lie algebra, which is a
linear Poisson structure.  In Refs. \cite{fl,ll} the case of a
compact semisimple group was considered and a a family of
 algebraic deformations of the coadjoint orbits was
constructed. This was done using the embedding of the coadjoint
orbit in the affine space as an affine algebraic, Poisson
subvariety.

\bigskip

In \cite{ko2} Kontsevich gives a  definition of {\it semiformal
deformations} of algebraic varieties. Basically, the deformed
algebra must have a filtration and at each step it should be
finitely generated.  Adapted to our notation the definition
 is the following.

\begin{definition}
Let $(\A,\{,\})$ be a finitely generated
Poisson algebra over a field $k$. We say that the
associative algebra $\Ah$ over $k[[h]]$ is a  semiformal
deformation of $\A$ if
\medskip

\noindent 1. $\A \simeq \Ah \otimes_{k[[h]]} k$ (i.e. $\Ah/(h)
\simeq \A$). We denote by $\pi:\Ah\rightarrow \Ah/(h) \cong \A$
the natural projection.

\noindent 2. There exists on $\Ah$  an exhaustive increasing
filtration, compatible with the product, and  admitting  a
splitting as a filtration of $k[[h]]$-modules. In other words
$\Ah=\cup_n \A^n$ with $\A^n \cdot \A^m \subset \A^{m+n}$, where
$\A^n$ are finitely generated free $k[[h]]$-modules, each a direct
addend of $\A_h$. This means that there exists a $k[[h]]$-module
$B^n$ such that $\A_h=A^n \oplus B^n$.

\noindent 3. $f_1 \star f_2-f_2 \star f_1=h\{\pi(f_1),\pi(f_2)\}\;
{\rmod}(h^2), \quad \forall f_1,f_2\in \Ah$.\hfill$\blacksquare$
\end{definition}

The concept of semiformal deformation differs from the formal one.
Let $\A$ be a commutative graded algebra and let
$$
\Ah=\A \otimes k[[h]], \hbox{ with } \A_n=\{f \in \A \otimes
k[[h]] \quad | \quad \deg f \leq n\}.
$$
Assume that we have a product in  $\Ah$ compatible with the
filtration. This is a prototype of a semiformal deformation. If
$\A$ is not finite dimensional the $k[[h]]$-module $\A \otimes
k[[h]]$ is strictly contained in $A[[h]]$, the $k[[h]]$-module
underlying a formal deformation of $\A$.

In some cases, given a formal deformation of $\A$, it is possible
to find a semiformal deformation sitting inside. This happens for
example when the the $k[[h]]$-module $\A \otimes k[[h]]$ is closed
under the star product.

In Ref.\cite{ko2}, \S 2.3,  Kontsevich discusses   this situation
for the algebra $\A=\Sym(V^*)$, where $V^*$ is a finite
dimensional vector space over $k$. It appears that only  Poisson
structures with degree up to 2 admit semiformal deformations.

Kontsevich also proposes in \S 4.1 the problem of finding a
semiformal deformation of regular coadjoint orbits of semisimple
Lie groups. Remember that in this case the Poisson structure is
linear. The problem was indeed already solved in \cite{fl,ll},
where formal deformations
 of the polynomial algebra  of coadjoint orbits were explicitly constructed.
 The deformations are also semiformal, in the sense that the subspace $\A
\otimes \C[[h]]$ is closed under the star product.

\bigskip

Another approach to the deformation of algebraic varieties is
taken in Ref. \cite{ye}. It is shown there that any smooth,
Poisson  algebraic variety (with some topological requirements)
admits a deformation. Also, such deformations are classified in
terms of formal Poisson structures (up to gauge equivalence), in
the same way that deformations of differential manifolds where
classified by Kontsevich \cite{ko1}. The basic idea is to endow
the smooth variety with a {\it differential trivialization} or
{\it \'etale coordinates}. In each open set the local result of
Ref. \cite{ko1} can be  reformulated in ring theoretic terms.
Then, one can glue the star products in different open sets using
a procedure analogous to the procedure that Fedosov used to show
the existence of star products on symplectic manifolds \cite{fe}.
This procedure was extended to general Poisson manifolds by
Cattaneo, Felder and Tomassini \cite{cft1} using the notions of
{\it formal geometry} by Gelfand and Kazhdan \cite{gk}. The
extension of these methods to the algebraic geometrical setting is
non trivial.

\bigskip

The work by Bezrukavnikov and Kaledin \cite{bk} deals also with
quantization in the algebro-geometric context. In particular, it
deals with symplectic smooth varieties, and shows that the Fedosov
quantization procedure can be translated into this context with
appropriate cohomological assumptions


\section{ \label{affinespace} Deformation quantization of affine space with a
Poisson structure}  We consider an open domain in $\R^n$ with an
arbitrary differential Poisson structure. Let us denote by
$\{x_i\}_{i=1}^n$ the coordinates in such open domain and let
$$\alpha=\alpha^{ij}(x)\partial_i\otimes\partial_j$$ denote the Poisson
structure. We want to briefly describe Kontsevich's local formula
for the star product canonically associated to $\alpha$ (for a
full description, we shall refer to the original paper
\cite{ko1}).

The star product is given in terms of certain admissible graphs,
each of which has associated a bidifferential operator
contributing to the sum (\ref{starproduct}) with an appropriate
weight.

The bidifferential operators at order $n$ are constructed with the
products of $n$ factors $\alpha$ acted by  partial derivatives
$\partial_i$ (up to 2n-2), and the indices contracted in
appropriate way. For example, at order $n=3$ one such operator
could be

$$\sum_{\scriptsize{\begin{matrix}i_1,i_2,i_3\\j_1,j_2,j_3\end{matrix}}}
(\partial_{i_2}\partial_{j_3}\alpha^{i_1j_1})(
\partial_{i_1}
\alpha^{i_2j_2})(\alpha^{i_3j_3})\partial_{i_3}\otimes\partial_{j_2}\partial_{j_1}$$

We are interested in Poisson structures on the whole affine space
$\R^n$ such that $\alpha^{ij}(x)$ are polynomial functions. Then,
the Poisson structure is algebraic and the operators $B_n$ in
(\ref{starproduct}) are bidifferential operators with polynomial
coefficients.  Moreover, denoting by $\deg(f)$ the maximum degree
of the polynomial $f$ we have
$\deg(B_m(f,g))\leq\deg(f)+\deg(g)+(p-2)n$, being $p=
\max_{(i,j)}(\deg(\alpha^{ij}))$. The star product of two
polynomials will be an infinite series in $h$ with coefficients in
$\R[x_1 \dots x_n]$, i.e. an element of $\R[x_1 \dots x_n][[h]]$.
Then, Kontsevich's star product is an algebraic star product.

In the following we will take the Poisson structure fixed and
Kontsevich's star product will be denoted simply by $\star$. We
will always work with the complexifications of the Poisson
structure and star product. In the affine space $\bA^n$ we choose
coordinates $\{x_i\}_{i=1}^n$ so $\C[\bA^n]=_{\defi} \C[x_1 \dots
x_n]$ (. The Poisson bracket is determined by its values on the
generators of this algebra,
$$\{x_i, x_j\}=\alpha^{ij}(x_1, \dots, x_n)$$
and we have
$$\star:\C[\bA^n][[h]]\times\C[\bA^n][[h]]\longrightarrow
\C[\bA^n][[h]].$$
\medskip

Our goal in this section is to give a presentation of the deformed
algebra $(\C[\bA^n],\star)$ in terms of generators and relations.
That is, we want to present it as a quotient of the ring of formal
power series in $h$ with coefficients in the full tensor algebra
$T(X_1 \dots X_n)$ generated by $X_1 \dots X_n$ and a two-sided
completed ideal $J$.
$$
(\C[\bA^n],\star)\simeq T(X_1 \dots X_n)[[h]]/J
$$

When the Poisson structure is of degree 0, ($p=0$, $\alpha_{ij} $
constant), then the star product is the Moyal star product.

 When the Poisson structure is homogeneous of
of degree 1  we have
$$\{x_i,x_j\}=c_{ij}^kx_k$$
and the star product algebra is isomorphic to the enveloping
algebra over $\C[[h]]$ \cite{ho,ko1} of the Lie algebra defined by
the structure constants $c_{ij}^k$.

The case of a homogeneous quadratic Poisson structure  is studied
by Kontsevich in Ref.\cite{ko2}, pg 11. In that paper it is shown
that the algebra $(\CAn \otimes \C[[h]], \star)$, which is
strictly smaller than the one we are considering,  is closed under
the star product and that it can be given in terms of generators
and (quadratic) relations. A similar presentation for Poisson
structures of higher degree is not provided in that paper.

\bigskip

Our strategy in solving this problem proceeds as follows:

\noindent 1. Let $I=(i_1,i_2,\dots, i_n)$ be a multiindex with
$i_j=1,\dots n$. We will prove that the ordered star monomials
i.e. the monomials
$$
x_{\star I}= x_{i_1} \star \dots \star x_{i_m} \qquad i_1 \leq
\dots \leq i_m
$$
are a basis for the $\C[[h]]$-module $\CAn[[h]]$. (Notice that
when we write $x_{i_1} \star \dots \star x_{i_m}$ omitting the
parenthesis we are making an implicit use of the associativity of
$\star$).

\noindent 2. Using part 1.  we will find  an algebra isomorphism
$$
\begin{CD}
\CAn[[h]] @>>>   T(X_1 \dots X_n)[[h]]/J \\
x_{i_1} \star \dots \star x_{i_m} @>>> X_{i_1} \dots X_{i_m}
\end{CD}
\qquad i_1 \leq \dots \leq i_m
$$
where $J$ is a two-sided completed ideal (completed in the
$h$-adic topology). $J$ is generated by the relations obtained
expressing the non ordered star monomials in terms of the ordered
ones.

\bigskip

We start applying the procedure to the $\C[[h]]$-module
$\CAn[[h]]/(h^N)\simeq\CAn[h]/(h^N)$. Then we will show that the
inverse limit of the algebras obtained is a formal deformation of
$\CAn$.

\begin{proposition}\label{basismodhn} Let $N\in \N$ be fixed. The ordered star monomials:
$$
x_{\star I}=x_{i_1} \star \dots \star x_{i_m}, \qquad I=(i_1 \dots
i_m), \qquad i_1 \leq \dots \leq i_m
$$
form a basis for the $\C[[h]]$-module $\CAn[h]/(h^N)$.
\end{proposition}

Proof. The ordered commutative monomials $x_I=x_{i_1} \dots
x_{i_m}$, $\;i_1 \leq \dots \leq i_m$ form a basis for $\CAn$. By
Kontsevich formula  we can express any ordered star monomial as an
infinite series in $h$ with coefficients in $\CAn$. Modulo $h^N$
the series becomes finite and can be rearranged as a finite
$\C[h]/(h^N)$-linear combination of the commutative monomials. We
denote it as
$$
x_{\star I}=\sum_J A_I^Jx_J,
$$
where $A_I^J\in \C[h]/(h^N)$.  Let us take an ordering in the set
of multiindices $I=(i_1 \dots i_m)$, $i_1 \leq \dots \leq i_m$.
With this ordering we denote by and $X_{\star}$ and  $X$ the
(infinite) column vectors of ordered star monomials and of
commutative monomials respectively. Then we have an infinite
$\C[h]/(h^N)$-linear system:
$$
X_{\star}=AX
$$
where $A$ is the infinite matrix with entries $A_I^J$. We notice
two crucial facts:

\medskip

\noindent (i) The matrix $A$ has only a finite number of entries
different from zero in each row (as we stated before, this is
because we are taking Kontsevich's formula modulo $h^N$).

\medskip

\noindent (ii) $A=\rId+hB$, since by Kontsevich formula the term
of order 0 in $h$ is the commutative product.

\medskip

\noindent This implies that $A$ is invertible. Its inverse in
$\CAn[[h]]/(h^N)$ can be written as
$$
A^{-1}=(\rId+hB)^{-1}=\sum_{m=0}^{N-1}(-1)^mh^mB^m
$$
Notice that $B^m$ makes sense because of property 1. In fact all
powers $B^m$ have only a finite number of entries different from
zero in each row. We have $X=A^{-1}X_{\star}$, hence the ordered
star monomials generate $\CAn[[h]]$.

\bigskip

We want now to prove linear independence. Let us assume that there
is a relation
$$
\sum_I a^I x_{\star I}=0 \mod(h^N),\quad \hbox{ i.e. } \quad\sum
a^I x_{\star I}=h^Nq$$ with $a^I\in \C[h]/(h^N)$ and  $q \in
\CAn[[h]]$.

 We have that  $\max_{\{I\}}(\deg(a^I)) \leq N$.
Specializing $h$ to zero we have
$$
\sum_I a^I(0) x_{I} =0 \quad \Rightarrow \quad a^I(0)=0
$$
Hence we have that our relation can be written:
$$
h \sum_I b^{I}  x_{\star I} = h^Nq
$$
with $a^I=hb^I$, with
$\max_{\{I\}}(\deg(b^I))<\max_{\{I\}}(\deg(a^I))$. Since
$\CAn[[h]]$ is an integrity domain this implies
$$
\sum b^{I}  x_{\star I} = h^{N-1}q.
$$
 We can again specialize to $h=0$, which will imply $ b^{I}(0)=0$, so
 $b^{I}=hc^{I}$, with $\max_{\{I\}}(\deg(c^I))<\max_{\{I\}}(\deg(b^I))$.
In each step we reduce the degree of the coefficients, so
repeating
 the argument a sufficient number of times we obtain our result.
 \hfill$\blacksquare$

\bigskip

Let $\cR \subset \CAn[h]/(h^N)$ be the (infinite) set of linear
relations expressing the non ordered star monomials in terms of
the ordered ones.
$$
x_{\star I}  =\sum_{J} d_I^J x_{\star J}, \quad j_1\leq\cdots \leq
j_m\;\;\mbox{and}\;\; i_1,\dots i_m \; \mbox{arbitrary}
$$

\medskip
Let  $T(X_1, \dots , X_n)[h]$ be the free tensor algebra generated
over $\C[h]$ by $X_1 \dots X_n$. We denote
$$
T_N=_{\defi} T(X_1, \dots , X_n)[h]/(h^N)
$$
Denote by $J_N$ the two-sided ideal generated in $T_N$ by the
relations in $\cR$ mod($h^N$), where we replace $\star$ with the
tensor multiplication.

\begin{proposition} \label{affinespaceprop} The $\C[h]/(h^N)$ linear morphism:
$$
\begin{CD}
  T_N/J_N @>\psi_N>> \CAn[h]/(h^N) \\
 X_{i_1}\otimes \cdots\otimes X_{i_r} @>>> x_{i_1} \star \dots \star x_{i_r} &
\end{CD}\qquad i_1 \leq \dots \leq i_r
$$
is well defined and it is an algebra isomorphism.
\end{proposition}

Proof. Consider the  surjective linear map
$$
\begin{CD}
 T_N @>\phi_N>>\CAn[h]/(h^N)  \\
 X_{i_1}\otimes \cdots\otimes X_{i_r} @>>> x_{i_1} \star \dots \star
x_{i_r}\qquad i_1,\dots i_r\quad \mathrm{arbitrary.}
\end{CD}
$$
 Notice that it
defines an algebra morphism. Moreover, $ J_N \subset \ker(\phi_N)
$. Hence the maps $\phi_N$ for all $N$ induce a family of
surjective algebra homomorphisms
$$
\begin{CD}
 T_N/J_N
@>\psi_N>> \C[x_1 \dots x_n][h]/(h^n) \\
 X_{i_1}\otimes \cdots\otimes X_{i_r} @>>>x_{i_1} \star \dots \star
x_{i_r}
\end{CD}$$
It is easy to see that ordered monomials in $T_N/J_N$ form a
basis. Because of the definition of $J_N$, it is clear that they
are a system of generators. Moreover, they are linearly
independent. In fact if there were a relation:
$$
\sum a_IX_I=0 \qquad \Rightarrow \qquad \sum a_I \psi_n(X_I) = 0
\qquad \Rightarrow \qquad a_I=0
$$
due to Proposition \ref{basismodhn}.

So we have obtained a surjective $\C[[h]]/(h^N)$-linear map
$\psi_N$ which  sends a basis into a basis. It is a linear
isomorphism. Since it also preserves the product, it  is also an
algebra isomorphism.\hfill$\blacksquare$

\bigskip

We want to consider now the limit $N\rightarrow \infty$.

\begin{theorem} \label{theoremaffinespace} Let $J_N$ be the family of ideals defined
above,
and let $ J=\underleftarrow{\lim} J_N $ be its inverse limit. Then
we have an algebra isomorphism
$$
(\CAn[[h]], \star) \cong T(X_1 ,\dots ,X_n)[[h]]/J.
$$
\end{theorem}

Proof. Consider the exact sequence:
$$  \begin{CD}
0 @>>> J_N @>>> T_N @>>> \CAn @>>> 0\end{CD}
$$
This is an exact sequence of inverse systems, i.e. it well behaves
with respect to the sequences defining the inverse systems.

In general an exact sequence of inverse systems does not
automatically give an exact sequence of the corresponding inverse
limits. However in this case it happens, since $J_n$ is a
surjective system i.e.
$$\begin{CD}
J_{N+1} @>>> J_{N}
\end{CD}$$
is surjective (see Ref. \cite{am} pg 104).

So we have:
$$\begin{CD}
0 @>>> \underleftarrow{\lim}  J_N @>>>
 \underleftarrow{\lim}  T_N @>>>
\underleftarrow{\lim}  \CAn @>>> 0\end{CD}
$$
which is what we wanted to prove.\hfill$\blacksquare$

In general, to obtain explicitly the generators of the ideal
$J=\underleftarrow{\lim}  J_n $ one needs to know the full star
product series. In the particular case of a quadratic Poisson
structure, they reduce to quadratic generators \cite{ko2}, but in
the more general case we don't even know if the number of
generators is finite.

Since $$\C[x_1,\dots x_n]\simeq T(X_1 ,\dots ,X_n)/(X_i\otimes
X_j-X_j\otimes X_i)$$ we have \begin{equation}\begin{CD}T(X_1
,\dots ,X_n)[[h]]/J@>>h\rightarrow 0>T(X_1 ,\dots
,X_n)/(X_i\otimes X_j-X_j\otimes
X_i)\label{limith0}\end{CD},\end{equation} which is another way of
expressing condition 2. in Definition \ref{formaldef}.


\section{ \label{algvariety} Deformation quantization of regular affine Poisson varieties}

We want to construct an explicit algebraic deformation
quantization of the ring of polynomial functions $\C[X]$ of a
regular affine Poisson variety. We assume that $X$ is a Poisson
subvariety of some Poisson structure defined in the  affine
ambient space  $\bA^n$. We denote by  $\I$ the ideal in
$\C[\bA^n]$ defining the variety $X$. By assumption, we have that
$\I$ is also a Poisson ideal, that is,
$$
\{\I,f\} \subset \I, \qquad \forall f \in \C[\bA^n].
$$
Let $\{p_1,\dots p_m\}$ be a basis of the ideal, so $\I=(p_1,\dots
p_m)$.  Let $\A_h=T(X_1 ,\dots ,X_n)[[h]]/J$ be the quantization
of the affine space as presented in Section \ref{affinespace}. We
want to construct an ideal $\I_h\in \A_h$ such that $\A_h/\I_h$ is
a deformation quantization of the Poisson algebra $\C[\bA^n]/\I$.

The general idea is the same than the one used for coadjoint
orbits of semisimple groups in Refs. \cite{fl,ll}. In that case we
had two advantages: one is the fact that the Poisson structure in
the ambient space is linear, which roughly allows to make
inductions on the degree of the polynomials; the other is that the
existence of the  action of the semisimple group provides some
tools that are not available in the general case. Nevertheless,
assuming some technical conditions it is possible to overcome the
difficulties.

The precise result is as follows:

\medskip

\begin{theorem} \label{bigtheorem}
Let $X$ be an affine Poisson variety with ideal $\I=(p_1 \dots
p_m)$ and Poisson structure induced from a Poisson structure in
$\C[\bA^n]$. Let $\A_h$ be an algebraic  deformation of the the
Poisson algebra $\C[\bA^n]$, as constructed in Section
\ref{affinespace}. Assume that

\noindent 1. The polynomials $\{p_1,\dots p_n\}$ are such that the
matrix $(dp_1,\dots dp_n)$ has maximal rank on the points of $X$.

\noindent 2. There exists liftings  $P_1,\dots ,P_m\in\A_h$  of
$p_1, \dots p_m$,
$$\begin{CD}P_1,\dots ,P_m@>>h\rightarrow 0>p_1, \dots p_m\end{CD}$$
 such that   the following left and right ideals
coincide:
$$
\I_h=(P_1, \dots ,P_2)_{\mathrm{left}}=(P_1, \dots,
P_m)_{\mathrm{right}}.
$$
Then $\A_h/ \I_h$ is an algebraic deformation quantization of
$\C[X]$.\hfill$\blacksquare$
\end{theorem}

Before going to the the proof, we want to  make remarks on the
hypothesis 2. For regular coadjoint orbits we can always fulfill
this condition. It is enough to take $p_i$ invariant, which is
always possible, and then to consider the Weyl map (or
symmetrizer) from polynomials into the enveloping algebra,
 $$W:\C[\bA^n]\simeq\mathrm{Sym}(X_1,\dots , X_n)\rightarrow U_h.$$
Then  $P_i=W(p_i)$ are in the center of the enveloping algebra and
become adequate liftings. For non regular orbits $p_1,\dots p_m$
can be chosen  spanning a finite dimensional representation of the
group, which is always possible. Then, lifting with the Weyl map
we obtain elements in the enveloping algebra satisfying condition
2. Even when condition 1. is not satisfied in this case, this
lifting was used in Ref. \cite{ll} to construct a deformation
quantization.

Another case where the lifting  is available is when the
generators $p_1,\dots, p_n$ are Casimirs of the Poisson structure,
that is
$$\{p_i,f\}=0,\qquad\forall f\in \C[\bA^n], \quad i=1,\dots n.$$
Indeed, it was found in Ref. \cite{cft1} (Theorem 5.1) the
explicit form of map $R:\C[x_1,\dots x_n][[h]]\rightarrow
\C[x_1,\dots x_n][[h]]$ which is the identity  when $h\rightarrow
0$ (a {\it quantization map}) that is an algebra isomorphism
between the Casimirs of the Poisson structure and the center of
the star product algebra. It is constructed in terms of the
$L_\infty$-morphism which gives the formality theorem \cite{ko1},
and when applied to polynomials gives a formal series in $h$ with
polynomial coefficients. Then it follows that $R(p_1), \dots ,
R(p_m)$ also satisfy condition 2. in Theorem \ref{bigtheorem}.

Note that R is an algebra isomorphism  but this is not necessary
to fulfill condition 2. The Weyl map, for example, is not an
algebra isomorphism.

Unfortunately, when the generators are not central, the same map
does not give an appropriate lifting. Nevertheless, it is very
likely that such lifting exists for every Poisson ideal, but  we
have not proved it in full generality.

Varieties defined by central elements are typically symplectic
leaves of the Poisson structure in the ambient space, or stacks of
such leaves. In particular, this extends our previous result on
coadjoint orbits of semisimple Lie groups to regular coadjoint
orbits of arbitrary (not necessarily semisimple) Lie groups.

\bigskip

We now return to the main result. We need only to prove that
$\A_h/\I_h$ is isomorphic, as a $\C[[h]]$-module to $\A/\I[[h]]$.
Then, properties 2. and 3. in Definition \ref{formaldef} are
immediate. We need some lemmas.

\begin{lemma}\label{olver}
Let $p_1, \dots p_m\in \C[x_1 \dots x_n]$ be such that the matrix
$(dp_1, \dots dp_m)$ has maximal rank on the points $p_1= \cdots
=p_m=0$. Then, if
$$
\sum_\alpha a_\alpha p_\alpha=0, \qquad a_\alpha\in  \C[x_1 \dots
x_n],
$$
 there exist elements $b_{\alpha\beta}\in  \C[x_1 \dots x_n]$ such that
$$
a_\alpha=\sum_j b_{\alpha\beta}p_\beta, \qquad\mathrm{with}\quad
b_{\alpha\beta}=-b_{\beta\alpha}
$$
\end{lemma}

Proof. This result can be found for $C^\infty$ functions in Ref.
\cite{ol}. We note that $b_ {ij}$ are not unique, because one can
always add a term $\tilde b_{ij}$ such that $\sum_j\tilde
b_{\alpha\beta}p_\beta=0$; that is, using the result in \cite{ol}
$$\tilde b_{\alpha\beta}=\sum_\gamma c_{\alpha\beta\gamma}p_\gamma, \qquad\mathrm{with}\quad
c_{\alpha\beta\gamma}=-c_{\alpha\gamma\beta}.
$$ It is
clear that if $p_\alpha$ and $a_\alpha$ are polynomials,
$b_{\alpha\beta}$ can also be chosen polynomial functions
\hfill$\blacksquare$
\bigskip

Next Lemma tells us that the $\C[[h]]$-module $\A_h/\I_h$ is
torsion free.

\begin{lemma} Let $X$, $p_1,\dots, p_m$ and $P_1,\dots, P_m$ be as in
 Theorem \ref{bigtheorem}. Let $\I_h=(P_1,\dots, P_m)\subset \A_h$.
Then, if  $hA\in \I_h$ also
 $A\in\I_h$.
\end{lemma}

Proof. Assume that  $hA \in \I_h$. Since $I_h$ is two-sided we can
write:
\begin{equation}
hA=\sum_\alpha A_\alpha P_\alpha. \label{torsion}
\end{equation}
Taking $h\rightarrow 0$ in this relation we get
$$\sum_ia_ip_i=0\quad\Rightarrow\quad a_\alpha=\sum_\beta b_{\alpha\beta}p_\beta,
 \quad \mathrm{with}\quad b_{\alpha\beta}=-b_{\beta\alpha}$$
by Lemma \ref{olver}. We can lift this relation,
$$\tilde A_\alpha=\sum_\beta B_{\alpha\beta }P_\beta, \quad \mathrm{with}\quad
B_{\alpha\beta}=-B_{\beta\alpha},$$ and it is clear that
$A_\alpha-\tilde A_\alpha=hC_\alpha$. By substituting in
(\ref{torsion}) we get
$$hA=\sum_\alpha (\sum_\beta B_{\alpha\beta}P_\beta+hC_\alpha)P_\alpha=h\sum_\alpha C_\alpha P_\alpha,$$
and so $A\in \I_h$.\hfill$\blacksquare$

\bigskip

We now want to define a set of elements $\B^*$ in $A_h$ whose
images in the quotient $A_h/I_h$ will turn out to be a topological
basis.

We consider a monomial basis for $\C[x_1 \dots x_m]$,
$$\{x_J=x_{j_1}\cdots x_{j_k}\},\qquad J=(j_1, \dots, j_k),\quad  j_1\leq \dots\leq j_k,
$$ and the associated topological basis of $\C[x_1 \dots
x_m][[h]]$ of star monomials
$$\{x_{\star J}=x_{j_1}\star\cdots\star x_{j_k}\},\qquad
J=(j_1, \dots, j_k),\quad  j_1\leq \dots\leq j_k, $$ (see Section
\ref{affinespace}).

 Let $B$ be  a set of multiindices with
the property the images by $\pi:\A\rightarrow \A/\I$ of the
elements in $\B=\{x_J\}_{J \in B}$  form a basis of $\A/\I$.
Consider the set
$$
\B^\star=\{x_{\star J}\}_{J \in B} \subset A_h.
$$
Let $\pi_h:\A_h \lra \A_h/\I_h$. We want to show that the elements
in $\pi_h(\B^\star)$ are a basis of $\A_h/\I_h$. We start by
showing that they are linearly independent.

\begin{lemma} \label{linind}
The image under the projection map $\pi_h:\A_h \lra \A_h/\I_h$ of
$\B^\star=\{x_{\star J}\}_{J \in B}$ is a linearly independent set
in $A_h/I_h$.

\end{lemma}

Proof. The argument of Proposition 3.11 in \cite{fl} works without
changes, we repeat it here for completeness.

 Suppose  that there exists  a linear relation among the elements
of $\B^\star$ and let  $G \in \I_h$  be such relation, $G=h^kF$
with $$\begin{CD}F@>>h\rightarrow 0>f\neq 0\end{CD}$$  for some
$k$. By Lemma \ref{torsion} we have that $F\in \I_h$, so
$$
F = \sum_\alpha A_\alpha P_\alpha ,\qquad \Rightarrow\qquad
\mathrm f=\sum_\alpha a_\alpha p_\alpha
$$
by taking $h\rightarrow 0$. But $f$ is a non trivial relation
between the ordinary monomials in $\B$, which is not possible
since $\B$ is assumed to be a basis of $\A/\I$. So we have proven
linear independence \hfill$\blacksquare$

\bigskip

To prove that the elements in Lemma $\pi_h(\B^\star)$ are a system
of generators we cannot use the same argument that appears in
Proposition 3.13 in Ref. \cite{fl}. The reason is that  having a
linear Poisson structure would allow us to  use an induction on
the degree of the polynomials,  so to  choose certain liftings  in
such way that the ``correction'' terms would have a degree
strictly smaller than the largest degree appearing in the order
$h^0$  term. Such induction argument is not possible for a Poisson
structure with arbitrary degree. Nevertheless, we have the
following

\begin{proposition}
Let the notation be as above. $\pi_h(\B^*)$ is a topological basis
for $\A_h/\I_h$.
\end{proposition}

Proof. We set
$$
(\A_h/\I_h)_N=(\A_h/\I_h)/(h^N).
$$
It is enough to show that $\pi(\B^*)$ is a basis for
$(\A_h/\I_h)_N$ for all $N$.

Let $P$ be a complement of $B$ in the set of multiindices, and let
$$\PP=\{x_J\}_{J \in P},\qquad \PP^\star=\{x_{\star J}\}_{J \in P}.$$
 We will prove that any monomial in $\pi(\PP^\star)$ is expressible in terms of
monomials  in $\pi_h(\B^\star)$.  This will clearly suffice since
we have proven that the monomials in $\pi_h(\B^\star)$ are
linearly independent (also in $(A_h/I_h)_N$).

We consider total lexicographic orderings in the sets $B$ and $P$.
Let $a\in \N$ be the position of certain  multiindex $K\in B$ and
let  $\mu \in \N$ be the position of certain  multiindex $J\in P$
with respect to the chosen orders. It will be convenient to denote
the respective monomials as
\begin{eqnarray*}
e_a=x_{K}, \; \in K\in \B, & v_\mu=x_{J}\; J \in P\\
e_{\star a}=x_{\star K}, \; \in K\in \B, & v_{\star \mu}=x_{\star
J}\; J \in P.
\end{eqnarray*}

Each  monomial in $\PP$ can be expressed as a linear combination
of monomials in $\B$ modulo an element in $\I$
$$v_\mu=\sum_ab_{\mu a}e_a+\sum_ic_{\mu i}(x)p_i,$$
where $b_{\mu a}\in \C$ and $c_{\mu i}(x)\in \C[\bA^n].$

We can lift  this relation to  the deformed algebra $\A_h$ in many
ways. Each star monomial $v_{\star\alpha} \in \pi(\PP^\star)$ can
be expressed as:
\begin{equation}
v_{\star\mu}=\sum_bB_{\mu a}(h)e_{\star a}+\sum_iC_{\mu
i}(x,h)\star P_i+h\sum A_{\mu\nu}(h)v_{\star\nu},
\label{classicalsystem}\end{equation} $B_{\mu a}(h),
A_{\mu\nu}(h)\in \C[[h]]$ and $C_{\alpha \mu}(x,h)\in
\C[\bA^n][[h]].$ Modulo $\I_h$ we have then a linear system
\begin{equation}\label{linearsystem}
\sum_{\nu}(\delta_{\mu\nu}-h
A_{\mu\nu}(h))v_{\star\nu}=\sum_aB_{\mu a}(h)e_{\star a},
\end{equation}
which in matrix form reads
$$
Dv_\star=Be_\star, \qquad D=\mathrm{Id} -hA.
$$
The infinite matrix $D$, modulo $h^N$,  has only a finite number
of entries non zero for each row, and modulo $h$ it is the
identity. Hence, by the same reasoning used in  Section
\ref{affinespace}, we can invert $D$, and its inverse in
$\CAn[h]/(h^N)$ can be written as
$$
D^{-1}=(\mathrm{Id}-hA)^{-1}=\sum_{m=0}^{N-1}h^mA^m.
$$
$D^{-1}$ has again a finite number of entries different from zero
in each row, so the multiplication $D^{-1}B$ is well defined and
$$
v_\star=D^{-1}Be_\star
$$
as we wanted to prove.

Making the inverse limit we have that $\pi(\B^\star)$ is a
topological basis of $\A_h/\I_h$. \hfill$\blacksquare$

\bigskip

This concludes the proof of Theorem \ref{bigtheorem}.

\section*{Acknowledgements}
R. F. wants to thank the Dipartimento di Fisica, Politecnico di
Torino and the Departament de F\'{\i}sica Te\`orica, Universitat
de Val\`encia for their kind hospitality during the realization of
this work. M. A. Ll. wants to thank the Dipartimento di
Matematiche, Universit\`a di Bologna and the Mathematics
Department of UCLA for their kind hospitality during the
realization of this work. V. S. V. wants to thank the Dipartimento
di Fisica, Politecnico di Torino and the INFN, Sezione di Torino
for their kind hospitality during the realization of this work.

The work of M. A. Ll. has been supported by the research grant BFM
2002-03681 from the Ministerio de Ciencia y Tecnolog\'{\i}a
(Spain) and from EU FEDER funds.

\end{document}